\documentclass[11pt]{amsart} 

\usepackage{tikz}
\usepackage{subfigure}
\usepackage{hyperref}

\newcommand\Cp{\mathbb{C}}
\newcommand\R{\mathbb{R}}
\newcommand\Z{\mathbb{Z}}
\newcommand\Id{\mathbf{Id}}
\newcommand\Span{\mathrm{Span}}
\newcommand\supp{\mathrm{supp}}
\newcommand\Col{\mathrm{Col}}
\newcommand\Row{\mathrm{Row}}
\newcommand\Prob{\mathbb{P}}
\newcommand\dist{\mathrm{dist}}
\newcommand\Sparse{\mathrm{Sparse}}
\newcommand\Comp{\mathrm{Comp}}
\newcommand\Incomp{\mathrm{Incomp}}
\newcommand\Exp{\mathbb{E}}
\newcommand\Var{\mathrm{Var}}
\newcommand\Net{\mathcal{N}}
\newcommand\cf{\mathcal{L}}
\newcommand\LCD{\mathrm{LCD}}
\newcommand\RLCD{\mathrm{RLCD}}
\newcommand\CLCD{\mathrm{CLCD}}
\newcommand\QCLCD{\mathrm{QCLCD}}
\newcommand\UD{\mathrm{UD}}

\newtheorem{problem}{Problem}

\numberwithin{equation}{section}

\begin{document}


\title{Quantitative invertibility of non-Hermitian random matrices}

\author{Konstantin Tikhomirov}

\address{School of Mathematics, 686 Cherry street, Atlanta GA 30332}
\email{ktikhomirov6@gatech.edu}
\thanks{This work was partially supported by the Sloan Research Fellowship and by NSF grant DMS 2054666}

\dedicatory{Dedicated to Prof. Nicole Tomczak-Jaegermann}

\begin{abstract}
The problem of estimating the smallest singular value of random square matrices
is important in connection with matrix computations and analysis of the spectral distribution.
In this survey, we consider recent developments in the study of quantitative invertibility
in the non-Hermitian setting, and review some applications of this line of research.
\end{abstract}

\maketitle


\section{Introduction}

Given an $N\times n$ ($N\geq n$) matrix $A$, its singular values are defined as square roots of the eigenvalues
of the positive semidefinite $n\times n$ matrix $A^* A$:
\[
s_i(A):=\sqrt{\lambda_i(A^* A)}, \quad i=1,2,\dots,n,
\]
where we assume the non-increasing ordering $\lambda_1(A^* A)\geq \lambda_2(A^* A)\geq\dots\geq \lambda_n(A^* A)$.
The classical Courant--Fischer--Weyl theorem provides a variational formula
\[
s_i(A)=\min\limits_{E:\;\dim(E)=n-i+1}\max\limits_{x\in E,\,\lVert x\rVert_2=1}\lVert Ax\rVert_2,\quad 1\leq i\leq n,
\]
where the minimum taken over all linear subspaces $E$ of the specified dimension.
In particular, \textit{the smallest} and \textit{the largest} singular values of $A$ can be computed as
\[
s_{\min}(A)=s_n(A)=\min\limits_{x:\,\lVert x\rVert_2=1}\lVert Ax\rVert_2,\quad s_{\max}(A)=s_1(A)=\max\limits_{x:\,\lVert x\rVert_2=1}\lVert Ax\rVert_2.
\]
Additionally, if the matrix $A$ is square ($N=n$) and invertible then $s_{\min}(A)=\frac{1}{s_{\max}(A^{-1})}$.

The magnitude of the smallest singular value of square random matrices has attracted much attention
due to the special role it plays in several questions of theoretical significance and in applications.
In particular, the ratio of the largest and smallest singular values of a square matrix ---
\textit{the condition number} --- is systematically used in numerical analysis as a measure of sensitivity
to round-off errors.
Further, for certain random matrix models, bounds on the spectral norm of the matrix' resolvent (or, equivalently, the smallest singular
value of diagonal shifts of the matrix) is a crucial point in the study of the spectral distribution.
We refer to Sections~\ref{s:num} and~\ref{s:spec} of the survey for a discussion of those directions.

In this survey, we consider \textit{quantitative invertibility}
of random \textit{non-Hermitian} square matrices, including matrices with independent entries and
adjacency matrices of random regular digraphs.
The main objective in that line of research is to obtain bounds on probabilities $\Prob\{s_{\min}(A)\leq t\}$ as a function of $t$,
of the dimension, and, possibly, of some parameters of the model under consideration,
such as the variance profile of the matrix or its mean.

One approach to the problem, which can be named analytical, is based on comparing the distribution of $s_{\min}(A)$
with the distribution of the smallest singular value of a corresponding Gaussian random matrix.
The latter is very well understood \cite{Edelman} 
since explicit formulas for the joint distribution of the singular values of Gaussian matrices
are available \cite{James}.
We refer to \cite{TV10a,CL2019} for results of that type.

Another approach, which is the focus of this survey,
falls in the category of \textit{non-asymptotic} methods \cite{RV survey} and
is based on a combination of
techniques originated within asymptotic geometric analysis. It often produces very strong
probability estimates although typically lacks the precision of the analytical methods.
The major features of that approach are (a) reducing the estimate for $s_{\min}$ to estimating distances between random vectors
and random linear subspaces associated with the matrix, and (b)
the use of concentration (Bernstein--type) and anti-concentration (Littlewood--Offord--type) inequalities.
Often, this approach also involves constructing discretizations of certain subsets of $\R^n$ or $\Cp^n$
(\textit{$\varepsilon$-nets})
and estimating their cardinalities.
We will give a description of the features by considering multiple examples
from the literature.

Because of some differences in methodology, and because we wish to emphasize the importance
of the matrix invertibility for numerical analysis and in the study of the spectral distribution,
this survey does not cover non-quantitative results on singularity of random matrices.
We note that estimating the singularity probability for several models of \textit{discrete} random matrices
is a major topic within the combinatorial random matrix theory \cite{K67,KKS1995,TV07a,BVW10,CTV,Hoi sym}.
In last few years there has been a significant progress in this research direction (also, as corollaries of quantitative results),
in particular, the problem of estimating the sigularity probability of adjacency matrices of random regular (di)graphs
\cite{JHuang,Meszaros,NM18}, of Bernoulli random matrices \cite{Tikh2020,LT2020,JSS2020a} and, more generally,
discrete matrices with i.i.d entries \cite{JSS2020b}, of random symmetric matrices \cite{CMMM,FJ19,CJMS,CJMSb}.
We refer to a recent survey \cite{Van survey new} for a discussion and further references.

The rest of the survey is organized as follows.
Sections~\ref{s:num} and~\ref{s:spec} provide motivation
for studying quantitative invertibility of non-Hermitian random matrices, and a brief account of known results.
In Section~\ref{s:method}, we give an overview the methodology, starting with the
result of Rudelson and Vershynin \cite{RV08a} as a main illustration.
We then discuss novel additions to the methodology made in the past 10 years, which allowed
to make progress on several important problems in the random matrix theory.
Finally, in Section~\ref{s:open}, we discuss some open problems.

Let us recall some notions which will be used further.

A random variable $X$ on $\R$ or $\Cp$ is called \textit{subgaussian} if $\Exp\,\exp(|X|^2/K^2)<\infty$
for some number $K>0$. The smallest value of $K$ such that $\Exp\,\exp(|X|^2/K^2)-1\leq 1$,
is called \textit{the subgaussian moment} of $X$.

Given a sequence of \textit{random} Borel probability measures $(\mu_m)_{m=1}^\infty$ and a random probability measure $\mu$ on $\Cp$,
we say that $\mu_m$ converge \textit{weakly in probability} to $\mu$ if for every bounded continuous function $f$ on $\Cp$,
\[
\lim_{m\to\infty}\Prob\bigg\{\Big|\int f\,d\mu_m-\int f\,d\mu\Big|>\varepsilon\bigg\}=0,\quad \forall\varepsilon>0.
\]

We will denote by $\lVert \cdot\rVert$ the spectral norm of a matrix. 
The standard Euclidean norm in $\R^n$ or $\Cp^n$ will be denoted by $\lVert \cdot\rVert_2$.
We will write $\dist(S,T)$ for the Euclidean distance between two subsets $S$ and $T$ of $\R^n$ or $\Cp^n$.
By $S^{n-1}(\R)$ or $S^{n-1}(\Cp)$ we denote the unit Euclidean sphere in $\R^n$ or $\Cp^n$, respectively.
The constants will be denoted by $C,c'$, etc.

\section{Quantitative invertibility in matrix computations}\label{s:num}

In this section, we discuss the importance of estimating
the smallest singular value in numerical analysis, and provide a brief overview of related results
on random matrices.

\subsection{The condition number in numerical analysis}
For an $n\times n$ invertible matrix $A\in \Cp^{n\times n}$, \textit{the condition number of $A$} is defined
as
\[
\kappa(A):=\lVert A\rVert\,\lVert A^{-1}\rVert=\frac{s_{\max}(A)}{s_{\min}(A)}.
\]
Consider a system of $n$ linear equations in $n$ variables, represented in the matrix-vector form as $Ax=b$.
If the system is \textit{well conditioned}, i.e the condition
number of the coefficient matrix $A$ is small,
a perturbation of the matrix or the coefficient vector does not
strongly affect the solution. In particular, the round-off errors in matrix computations
such as the Gaussian elimination, do not significantly distort the solution vector.

As an example of well known theoretical guarantees, we mention
an estimate on \textit{the relative distance} between the solution of
$Ax=b$ and the solution of a perturbed system
\[
(A+F)y=(b+f).
\]
The terms $F\in \Cp^{n\times n}$ and $f\in\Cp^n$ can be thought of as consequences of
measurement or round-off errors.
It is not difficult to check
that under the assumption that $\delta:=\max(\frac{\lVert F\rVert}{\lVert A\rVert},\frac{\lVert f\rVert_2}{\lVert b\rVert_2})$
is small,
the relative distance $\frac{\lVert y-x\rVert_2}{\lVert x\rVert_2}$ satisfies
\[
\frac{\lVert y-x\rVert_2}{\lVert x\rVert_2}=O\big(\delta\,\kappa(A)\big)
\]
(see, in particular, \cite[Section~2.6.2]{MatrixComp4},
\cite[Section~4]{Smale1985}).
In the specific setting when the system $Ax=b$
is solved using the Gaussian elimination with partial pivoting and the perturbation of the system
is due to round-off errors,
Wilkinson \cite{Wilkinson} showed that the relative distance between the computed and the actual solution can be bounded above by
\[
n^{O(1)}\,\varepsilon\,\kappa(A)\,\rho.
\]
Here \textit{the growth factor} $\rho$ is defined as $\rho:=\frac{\max_{k=0,1,\dots;\,i,j\leq n}|a_{ij}^{(k)}|}
{\max_{i,j\leq n}|a_{ij}^{(k)}|}$,
with $a_{ij}^{(k)}$ being the $(i,j)$--th element of the matrix $A^{(k)}$ obtained from $A$ after $k$ iterations
of the Gaussian elimination process, and $\varepsilon$ is the precision of the machine (see also \cite{TS1990,Sankar}).

\medskip

Whereas the condition number of $A$ characterizes sensitivity of the corresponding SLE to small perturbations,
\textit{the eigenvector condition number} quantifies stability of the spectrum and eigenvectors of $A$.
The eigenvector condition number of a diagonalizable matrix $A\in \Cp^{n\times n}$ is defined as
\[
\kappa_V(A):=\min_{W\in\Cp^{n\times n}:\,W^{-1}A W\mbox{ is diagonal}}\kappa(W)=
\min_{W\in\Cp^{n\times n}:\,W^{-1}A W\mbox{ is diagonal}}\frac{s_{\max}(W)}{s_{\min}(W)}.
\]
Clearly, $\kappa_V(A)=1$ if and only if $A$ is \textit{unitarily diagonalizable} (normal).
A classical stability result for a matrix spectrum using the eigenvector condition number
is the Bauer--Fike theorem \cite{BF1960}.
According to the theorem, given a diagonalizable matrix $A$ and its perturbation $A+F$,
the distance between any eigenvalue $\mu$ of $A+F$ and the spectrum of $A$ can be estimated as
\[
\min_{\lambda\in \textrm{Spec}(A)}|\mu-\lambda|\leq \kappa_V(A)\lVert F\rVert.
\]
Moreover, stability of matrix functions under pertubations of the argument
can be quantified using the eigenvector condition number (see \cite[Section~3.3]{Higham2008}).
Here, we refer to a related line of research dealing with \textit{the approximate diagonalization}
of matrices --- approximating a matrix with one with a small eigenvector condition number
(see \cite{Davies,BKMS2019,BGKS2020,JSS2020} and references therein).
A connection between $\kappa_V(A)$ and quantitative invertibility of diagonal shifts of $A$
is established through the notion of a pseudospectrum.
\textit{An $\varepsilon$--pseudospectrum of $A$} --- $\textrm{Spec}_{\varepsilon}(A)$ --- is defined as the set of all points $z\in\Cp$
with $s_{\min}(A-z\,\Id)<\varepsilon$. It can be shown (see \cite[Lemma~9.2.11]{DaviesBook})
that for a diagonalizable matrix $A$ with $D$ being a corresponding diagonal matrix,
$
\textrm{Spec}(D)+\kappa_V(A)^{-1}\,\varepsilon\,U\subset
\textrm{Spec}_{\varepsilon}(A)\subset \textrm{Spec}(D)+\kappa_V(A)\,\varepsilon\,U$, where $U$ is the unit disk of the complex plane.

\subsection{Related results on random matrices}
Randomness is a natural approach to simulate typical matrices observed in applications.
For example, the LINPACK benchmark for measuring the computing power involves systems of linear equations
with a randomly generated coefficient matrix \cite{LINPACK}.
Condition numbers of random square matrices with the computational perspective
were first considered by von Neumann and Goldstine \cite{vNG}.
Rigorous results were obtained much later, notably by Edelman \cite{Edelman}
for Gaussian random matrices  (see also Szarek \cite{Szarek}).
We note here that for \textit{sufficiently dense} random matrices with i.i.d entries satisfying certain moment conditions,
estimating the largest singular values up to a constant multiple can be accomplished by a simple combination
of Bernstein--type inequalities and an $\varepsilon$--net argument (see, for example \cite{RV survey}),
and with precision up to $(1\pm o(1))$ multiple via the trace method \cite{Geman,YBK,Seginer}.
Further, we only discuss estimates for the smallest singular value.

The average-case quantitative analysis of the matrix invertibility, when
a typical matrix is modeled as a random matrix with independent entries with matching first two moments,
has been developed in multiple works. We refer, in particular, to papers \cite{TV10a,CL2019}
employing the analytical approach, as well as works
\cite{Rud2008,TV2009,RV08a,RV08b,RT15,BR16,BR threshold,
Luh,Tikh shifted,Livshyts,Tikh2020,LTV,Tatarko,LT2020,JSS2020a,JSS2020b,Huang}
based on reduction to distance estimates
and on the use of concentration/anti-concentration inequalities. Some of those results are mentioned below.

In \cite{RV08a}, Rudelson and Vershynin showed that given a random $n\times n$ matrix $A$
with i.i.d real entries of zero mean, unit variance and a bounded \textit{subgaussian moment},
the smallest singular value of $A$ satisfies
\[
\Prob\{s_{\min}(A)\leq n^{-1/2}\,t\}\leq C(t+c^n),\quad t>0,
\]
where the constants $C>0$ and $c\in(0,1)$ may only depend on the sugaussian moment
(in fact, the statement is preserved if $A$ is shifted by a non-random matrix with the spectral norm of order $O(\sqrt{n})$).
The moment assumptions and the requirement that the entries are equidistributed were relaxed in later works \cite{RT15,Livshyts,LTV}.
On the other hand, in the special case of a matrix $A$ with i.i.d entries taking values $+1$ and $-1$ with probability $1/2$,
it was proved in \cite{Tikh2020} that for any $\varepsilon>0$,
\[
\Prob\big\{s_{\min}(A)\leq n^{-1/2}\,t\big\}\leq Ct+C(1/2+\varepsilon)^n,\quad t>0,
\]
where $C>0$ is only allowed to depend on $\varepsilon$ (see the introduction to \cite{Tikh2020}
as well as \cite{Van survey new} for a discussion of this result in the context of the combinatorial random matrix theory).
Yet stronger result is available when $A$ has i.i.d \textit{discrete} entries which are not uniformly distributed on their support
\cite{JSS2020b}: for every $\varepsilon>0$ and assuming $n$ is sufficiently large,
\[
\Prob\big\{s_{\min}(A)\leq n^{-1/2}\,t\big\}\leq Ct+(1+\varepsilon)
\Prob\big\{\mbox{Two rows or columns of $A$ are colinear}\big\},\quad t>0,
\]
with $C>0$ depending only on the individual entry's distribution (see \cite{JSS2020b} for the statement in its full strength).
In the setting when $A$ has i.i.d Bernoulli($p$) entries and $p$ is allowed to depend on $n$,
its was shown in \cite{BR threshold,LT2020,Huang} that, as long as $p\leq c$
for a small universal constant $c>0$, for every $\varepsilon>0$ and assuming $n$ is sufficiently large,
\[
\Prob\big\{s_{\min}(A)\leq n^{-C}\,t\big\}\leq t+(1+\varepsilon)
\Prob\big\{\mbox{A row or a column of $A$ is zero}\big\},\quad t>0,
\]
where $C>0$ is a universal constant. We refer to \cite{Huang} for a generalization to matrix rank estimates, as well as
work \cite{JSS2020a} for sharp bounds in the setting of constant $p\in(0,1/2)$ and \cite{BR16}
for stronger quantitative estimates in a certain range for the parameter $p$.

\medskip

Put forward by Spielman and Teng \cite{ST02},
\textit{the smoothed analysis} of the condition number is concerned with quantitative invertibility of a typical matrix
in a small neighborhood of a fixed matrix (with possibly a very large spectral norm).
A basic probabilistic model of that type is of the form $A+M$, where $M$ is a non-random matrix, and $A$ has i.i.d entries.
The result of Sankar--Spielman--Teng \cite{SST06}
provided a small ball probability bound for the smallest singular value of a shifted \textit{Gaussian} real random matrix
with i.i.d standard normal entries,
\textit{independent} from the shift $M$:
\[
\Prob\big\{s_{\min}(A+M)\leq t\,n^{-1/2}\big\}\leq C\,t,\quad t>0,
\]
for a certain universal constant $C>0$ (see also \cite[Section~2.3]{BKMS2019}).
Analogous estimates for a broader class of random matrices with
continuous distribution was later obtained in \cite{Tikh shifted}.
On the other hand,
it was observed that for certain discrete random matrices, such as random sign (Bernoulli) matrices, no shift-independent
small ball probability bounds for $s_{\min}(A+M)$ are possible \cite{TV10c,Tikh shifted,JSS2020 smooth}.
In particular, it is shown in \cite{JSS2020 smooth} that, assuming $A$ has i.i.d entries taking values $\pm 1$ with probability $1/4$
and zero with probability $1/2$, every $L\geq 1$, and every positive integer $K$,
\[
\sup_{M:\,\lVert M\rVert \leq n^{L}}\Prob\big\{s_{\min}(A+M)\leq C\, n^{-KL}\big\}\geq c\,n^{-K(K-1)/4},
\]
where $C,c>0$ may only depend on $L$ and $K$.
The smoothed analysis of the matrix condition number for discrete distributions was carried
out in works \cite{TV STOC, TV10c,Jain2021a,JSS2020 smooth} (see also references therein).
The following result was proved in \cite{TV10c}.
Let $K,B,\varepsilon>0$ and $L\geq 1/2$ be arbitrary parameters.
Then, for all sufficiently large $n$, given an $n\times n$ random matrix $A$ with i.i.d centered entries
of unit variance and the subgaussian moment bounded above by $B$,
and given a non-random matrix $M$ with $\lVert M\rVert\leq n^L$,
one has
\[
\Prob\big\{s_{\min}(A+M)\leq n^{-(2K+1)L}\big\}\leq n^{-K+\varepsilon}.
\]
In \cite{JSS2020 smooth}, it is shown that the above small ball probability bound can be significantly improved to match the
average-case result of Rudelson and Vershynin \cite{RV08a}, under the assumption that a positive fraction of the
singular values of $M$ are of order $O(\sqrt{n})$.
More specifically, for every $\tilde c\in(0,1)$ and $\tilde C>0$, and any fixed matrix $M$ with $s_{n-\lfloor \tilde cn\rfloor}(M)\leq \tilde C\sqrt{n}$,
one has
\[
\Prob\big\{s_{\min}(A+M)\leq t\,n^{-1/2}\big\}\leq C(t+c^{n}),
\]
where $C>0$ and $c\in(0,1)$ may only depend on $\tilde c$, $\tilde C$, and the subgaussian moment $B$.
Under much weaker assumptions on the shift $M$, though at a price of
precision, quantitative bounds for $s_{\min}(A+M)$ were obtained in \cite{Jain2021a}.

\section{Invertibility and spectrum}\label{s:spec}

Given a square $n\times n$ matrix $A_n$, denote by $\mu_{A_n}$ its normalized spectral measure (spectral distribution):
$$
\mu_{A_n}:=\frac{1}{n}\sum_{i=1}^n \delta_{\lambda_i(A_n)}.
$$
For real and complex Gaussian matrices with i.i.d standard entries (\textit{the Ginibre ensemble})
explicit formulas for the joint distribution of the eigenvalues are known
\cite{Ginibre,EdelmanJMA,LS1991}.
Those, in turn, were used by Mehta \cite{Mehta}, Silverstein (unpublished; see \cite[Section~3]{BC2011})
and Edelman \cite{EdelmanJMA} to derive convergence results for the spectral distribution in the Gaussian case.

In the non-Gaussian setting, when no similar formulas are available,
Girko \cite{Girko} proposed a \textit{Hermitization} argument based on the identity
\begin{align*}
\frac{1}{n}\sum_{i=1}^n\log|z-\lambda_i(A_n)|
=\frac{1}{n}\log\sqrt{\det((A_n-z\,\Id)(A_n-z\,\Id)^*)}=\frac{1}{n}\sum_{i=1}^n\log s_i(A_n-z\,\Id),
\end{align*}
which relates the spectrum to the singular values of the matrix resolvent.
A modern form of the argument can be summarized as follows (see \cite[Lemma~4.3]{BC2011}
as well as \textit{the replacement principle} in \cite{TV circ}).
Assume that a sequence of random matrices $(A_n)_{n=1}^\infty$ is such that for almost every $z\in\Cp$,
the sequence of measures $\mu_{\sqrt{(A_n-z\,\Id)(A_n-z\,\Id)^*}}$ converges weakly in probability
to a non-random probability measure $\mu_z$. Assume further that the logarithm is \textit{uniformly integrable in probability}
with respect to $(\mu_{\sqrt{(A_n-z\,\Id)(A_n-z\,\Id)^*}})_{n=1}^\infty$ for almost every $z\in\Cp$, that is
\begin{equation}\label{1-9481-498-}
\lim_{t\to\infty}\sup_n \Prob\bigg\{\frac{1}{n}\sum_{i\leq n:\,|\log s_i(A_n-z\,\Id)|>t}|\log s_i(A_n-z\,\Id)|>\varepsilon\bigg\}=0,\quad
\forall \varepsilon>0.
\end{equation}
Then there is a measure $\mu$ on $\Cp$ such that the sequence $(\mu_{A_n})_{n=1}^\infty$
converges to $\mu$ weakly in probability; moreover, the measure $\mu$ can be characterized in terms of $(\mu_z)_{z\in\Cp}$.
We refer to \cite{BC2011} for proofs as well as a detailed historical account of the study of the spectral distribution
of non-Hermitian random matrices, up to 2000-s.

In view of the uniform integrability requirement \eqref{1-9481-498-}, strong quantitative estimates for small singular values of
matrices $A_n-z\,\Id$ are an essential part of the Hermitization argument.
In the setting of matrices with i.i.d non-Gaussian entries,
first rigorous estimates on the small singular values of $A_n-z\,\Id$ sufficient for the argument to go through were obtained
for a class of continuous distributions by Bai \cite{Bai circ}, who applied the estimates
to study the limiting spectral distribution in that setting. As the techniques to quantify invertibility of 
more general classes of matrices
became available through the works of Tao--Vu \cite{TV2009}, Rudelson \cite{Rud2008}, and Rudelson--Vershynin \cite{RV08a},
the result of Bai was consecutively generalized in works \cite{GT circ,PZ circ,TV circ -,TV circ}.
The \textit{strong circular law} under minimal moment assumptions proved in \cite{TV circ} can be formulated as follows.
Let $\xi$ be a complex valued random variable of zero mean and unit absolute second moment,
and let $(A_n)_{n=1}^\infty$ be a sequence of random matrices, where each $A_n$ is $n\times n$
with i.i.d entries equidistributed with $\xi$.
Then the sequence of spectral distributions $(\mu_{\frac{1}{\sqrt{n}}A_n})_{n=1}^\infty$ converges weakly almost surely
to the uniform probability measure on the unit disk of the complex plane.

\smallskip

In the context of the circular law, a most studied model of \textit{sparse} random matrices is of the form $A_n=B_n \odot M_n$,
where $B_n$ is the random matrix with i.i.d Bernoulli($p_n$) entries, $M_n$ is independent from $B_n$
and has i.i.d entries equidistributed with a random variable
$\xi$ of unit variance, and ``$\odot$'' denotes the Hadamard (entry-wise) product of matrices.
In the regime $p_n\geq n^{-1+\varepsilon}$ for a fixed $\varepsilon>0$, the (weak) circular law
has been established in \cite{Wood circ} following earlier works \cite{GT circ,TV circ -} dealing with additional moment assumptions.

In yet sparser regime, estimating the smallest singular value of $A_n-z\,\Id$ presents significant challenges,
and further progress has only been made recently in \cite{BR circ,RT circ}.
In \cite{RT circ}, it is proved that, assuming $\xi$ is a real valued random variable of unit variance, $n\,p_n\leq n^{1/8}$, and
$n p_n$ tends to infinity with $n$, and assuming the matrices $A_n=B_n \odot M_n$ are defined as in the previous paragraph,
the sequence of spectral distributions $\big(\mu_{\frac{1}{\sqrt{n p_n}}A_n}\big)_{n=1}^\infty$
converges weakly in probability to the uniform measure on the unit disk of $\Cp$.
A central technical result of \cite{RT circ}
is the following quantitative bound for $s_{\min}(A_n-z\,\Id)$: under the assumption that $|z|\leq np_n$ and $|\Im(z)|\geq 1$,
\[
\Prob\big\{s_{\min}(A_n-z\,\Id)\leq \exp(-C\log^3 n)\big\}\leq C\,(np_n)^{-c},
\]
where $C,c>0$ may only depend on the c.d.f of $\xi$.

Quantitative invertibility and spectrum of adjacency matrices of random regular directed graphs have
been considered in multiple works in past years \cite{C14,Cook circ random,Cook circ,BCZ circ,LLTTY sing,LLTTY trans,LLTTY circ}.
Given integers $n$ and $d$, a \textit{$d$--regular digraph} on vertices $\{1,2,\dots,n\}$
is a directed graph in which every vertex has $d$ incoming edges and $d$ outgoing edges.
Here, we focus on the model when no multiedges are allowed but the graph may have loops
(the latter condition is not conventional).
For each $n$, denote by $A_{n,d}$ the adjacency matrix of a random graph
uniformly distributed on the set of all $d$--regular digraphs on $\{1,2,\dots,n\}$ (we allow $d$ to depend on $n$).
First results on invertibility for this model were obtained by Cook \cite{C14}.
The circular law for the sequence of spectral measures $\big(\mu_{\frac{1}{\sqrt{d(1-d/n)}}A_{n,d}}\big)_{n=1}^\infty$
has been established in \cite{Cook circ} under the assumption $\min(d,n-d)\geq \log^{96} n$.
Later, in works \cite{LLTTY sing,LLTTY trans,LLTTY circ}, the range $\omega(1)=d\leq \log^{96} n$ was treated.
Either of the two results relies heavily on estimates of the smallest singular values of $A_{n,d}-z\,\Id$.
In particular, the main theorem of \cite{LLTTY sing} is the following statement:
assuming $C\leq d\leq n/\log^2 n$ and $|z|\leq d/6$,
\[
\Prob\big\{s_{\min}(A_{n,d}-z\,\Id)< n^{-6}\big\}\leq \frac{C\log^2 d}{\sqrt{d}},
\]
where $C>0$ is a universal constant.

\smallskip

Invertibility of \textit{structured} random matrices and applications to the study of limiting spectral distribution
have been considered, in particular, in \cite{RZ2016,C16,Cooketal1,Cooketal2,JJLR}.
A basic model of interest here is of the form $A_n=U_n\odot M_n-z\,\Id$, where $M_n$ is a matrix
with i.i.d entries of zero mean and unit variance, $z\in\Cp$ is some complex number, $U_n$
is a non-random matrix with non-negative real entries
encoding \textit{the standard deviation profile}, and ``$\odot$'' denotes the Hadamard (entry-wise)
product of matrices. Note that $A_n$ has mutually independent entries, with $\sqrt{\Var\, a_{ij}}=u_{ij}$, $1\leq i,j\leq n$.
In \cite{RZ2016}, invertibility (and, more generally, the singular spectrum) of $U_n\odot M_n$
was studied in connection with the problem of estimation of matrix permanents.
In particular, strong quantitative bounds on $s_{\min}(U_n\odot M_n)$ were obtained in the setting when $M_n$
is the standard real Gaussian matrix, and $U_n$ is a \textit{broadly connected} profile (see \cite[Section~2]{RZ2016}).
A significant progress in the study of structured random matrices was made by Cook in \cite{C16},
who extended the result of \cite{RZ2016} to non-Gaussian matrices, and obtained a polynomial lower bound on
$s_{\min}(U_n\odot M_n-z\,\Id)$ under very general assumptions on $U_n$. Namely, assuming that all entries of $U_n$
are in the interval $[0,C]$, that $z\in [c\sqrt{n},C\sqrt{n}]$ for some constants $c,C>0$, and that the entries of $M_n$
have a bounded $(4+\varepsilon)$--moment, the main result of \cite{C16} asserts that
\[
\Prob\big\{s_{\min}(U_n\odot M_n-z\,\Id)\leq n^{-\beta}\big\}\leq n^{-\alpha}
\]
for some $\alpha,\beta>0$ depending only on $c,C,\varepsilon$, and the value of the $(4+\varepsilon)$--moment.
In \cite{Cooketal1}, this estimate was applied to derive limiting laws for the spectral distributions, under some
additional assumptions on $U_n$. One of the results of \cite{Cooketal1} is \textit{the circular law for doubly
stochastic variance profiles}: provided that
$\sum_{i=1}^n (U_n)^2_{ij}=\sum_{i=1}^n (U_n)^2_{ji}=n$, $1\leq j\leq n$,
and $\sup_n\max_{i,j}(U_n)_{ij}<\infty$, the sequence of spectral distributions $\big(\mu_{\frac{1}{\sqrt{n}}U_n\odot M_n}\big)_{n=1}^\infty$
converge weakly in probability to the uniform measure on the unit disc of $\Cp$.

The setting of \textit{sparse} structured matrices is not well understood. 
For results in that direction, we refer to a recent paper \cite{JJLR} dealing with
invertibility and spectrum of \textit{block band matrices}.

\section{Methodology}\label{s:method}

We start this section with a brief outline of \cite{RV08a}
which will serve as a canonical illustration of non-asymptotic methods.
The proof of the main theorem in \cite{RV08a} relies on four major components: sphere partitioning, invertibility via distance,
$\varepsilon$--net arguments, and Littlewood--Offord--type inequalities.

Let $A$ be an $n\times n$ matrix with i.i.d real entries of zero mean and unit variance,
and assume for simplicity that the entries are $K$--subgaussian for some constant $K>0$.
We recall that a vector $x\in\R^n$ is called \textit{$m$--sparse} if the size of its support is at most $m$.
We will denote the set of all $m$--sparse vectors by $\Sparse_n(m)$.
The proof of \cite[Theorem~3.1]{RV08a} starts with splitting $S^{n-1}(\R)$
into sets of \textit{compressible} and \textit{incompressible} vectors,
\begin{align*}
\Comp_n(\delta,\rho)&:=\big\{x\in S^{n-1}(\R):\;\dist(x,\Sparse_n(\delta n))< \rho\big\};\\
\Incomp_n(\delta,\rho)&:=\big\{x\in S^{n-1}(\R):\;\dist(x,\Sparse_n(\delta n))\geq \rho\big\}.
\end{align*}
Here, $\delta,\rho\in(0,1)$ are small constants.
The variational formula for $s_{\min}(A)$ allows to write
\begin{align*}
\Prob\big\{s_{\min}(A)\leq s\big\}\leq
&\Prob\big\{\lVert Ax\rVert_2\leq s\mbox{ for some $x\in \Comp_n(\delta,\rho)$}\big\}\\
+
&\Prob\big\{\lVert Ax\rVert_2\leq s\mbox{ for some $x\in \Incomp_n(\delta,\rho)$}\big\},\quad s>0.
\end{align*}
If both $\delta$ and $\rho$ are sufficiently small, the set of compressible vectors has small \textit{covering numbers},
which allows to apply an $\varepsilon$--net argument. More specifically, it can be checked that for every $\varepsilon\in(3\rho,1/2]$
there is a discrete subset $\Net\subset \Comp_n(\delta,\rho)$ of size at most $\big(\frac{C}{\varepsilon\delta}\big)^{\delta n}$
such that for every $x\in \Comp_n(\delta,\rho)$, we have $\dist(x,\Net)\leq\varepsilon$ (i.e, $\Net$ is
an $\varepsilon$--net in $\Comp_n(\delta,\rho)$ with respect to the Euclidean metric).
Consequently, for every $L>0$,
\begin{align*}
\Prob&\big\{\lVert Ax\rVert_2\leq s\mbox{ for some $x\in \Comp_n(\delta,\rho)$}\big\}\\
&\leq \Prob\big\{\lVert Ay\rVert_2\leq s+\varepsilon\,L\sqrt{n}\mbox{ for some $y\in \Net$}\big\}
+\Prob\big\{\lVert A\rVert> L\sqrt{n}\big\}\\
&\leq |\Net|\,
\sup\limits_{z\in S^{n-1}(\R)}\Prob\big\{\lVert Az\rVert_2\leq s+\varepsilon\,L\sqrt{n}\big\}
+\Prob\big\{\lVert A\rVert> L\sqrt{n}\big\}.
\end{align*}
For any $z\in S^{n-1}(\R)$, the vector $Az$ has i.i.d subgaussian components of unit variances,
and a standard Laplace transform argument implies that, as long as $s+\varepsilon\,L\sqrt{n}$ is much less than $\sqrt{n}$,
the probability $\Prob\big\{\lVert Az\rVert_2\leq s+\varepsilon\,L\sqrt{n}\big\}$ is exponentially small in $n$.
Moreover, for a sufficiently large constant $L$, the probability $\Prob\big\{\lVert A\rVert> L\sqrt{n}\big\}$ 
is exponentially small in $n$. Therefore, an appropriate choice of parameters $\delta,\rho,\varepsilon,L$ yields
\[
\Prob\big\{\lVert Ax\rVert_2\leq s\mbox{ for some $x\in \Comp_n(\delta,\rho)$}\big\}\leq 2\exp(-cn),\quad s=o(\sqrt{n}).
\]
We refer to \cite{RV08a} as well as \cite{RV survey} for details regarding the above computations.
Let us note also that the idea of sphere partitioning was applied a few years earlier in paper \cite{LPRT}
dealing with rectangular random matrices.

The incompressible vectors are treated using the \textit{invertibility via distance} argument,
which is based on the observation that for any incompressible vector $x$, a constant proportion of its components
are of order $\Omega(n^{-1/2})$ by the absolute value.
For every $1\leq i\leq n$, denote by $H_i(A)$
the linear span of columns of $A$ except the $i$-th:
\[
H_i(A):=\Span\{\Col_j(A),\;j\neq i\}.
\]
Then for arbitrary vector $x$ and arbitrary ``threshold'' $\tau>0$ with $\{i:\,|x_i|\geq \tau\}\neq\emptyset$ we have
\[
\lVert Ax\rVert_2\geq \max_{1\leq i\leq n}\big(|x_i|\,\dist(\Col_i(A),H_i(A))\big)
\geq \tau\max_{i:\,|x_i|\geq \tau}\dist(\Col_i(A),H_i(A)),
\]
and hence for any $s>0$, $
\textbf{1}_{\{\lVert Ax\rVert_2\leq s\}}\leq \frac{1}{|\{i:\,|x_i|\geq \tau\}|}\sum_{i=1}^n\textbf{1}_{\{\dist(\Col_i(A),H_i(A))\leq s/\tau\}}$.
This, combined with Markov's inequality and the fact that every $(\delta,\rho)$--incompressible vector is \textit{spread} ---
has at least $\delta n$ components of magnitude at least $\rho n^{-1/2}$ --- gives for $t>0$,
\begin{equation}\label{310984-198}
\Prob\big\{\exists\,x\in \Incomp_n(\delta,\rho):\; \lVert Ax\rVert_2\leq t\,n^{-1/2}\big\}
\leq \frac{1}{\delta n}\sum_{i=1}^n\Prob\{\dist(\Col_i(A),H_i(A))\leq t/\rho\}
\end{equation}
(see \cite[Lemma~3.5]{RV08a}).
Since the distribution of $A$ is invariant under columns permutations, the last relation can be rewritten as
\begin{align*}
\Prob\big\{\lVert Ax\rVert_2\leq t\,n^{-1/2}\mbox{ for some $x\in \Incomp_n(\delta,\rho)$}\big\}
&\leq \frac{1}{\delta}\Prob\{\dist(\Col_n(A),H_n(A))\leq t/\rho\}\\
&\leq\frac{1}{\delta}\Prob\{|\langle \Col_n(A),Y_n(A)\rangle|\leq t/\rho\},
\end{align*}
where $Y_n(A)$ denotes a unit normal to $H_n(A)$ measurable with respect to $\sigma(H_n(A))$.

The most involved part of \cite{RV08a} is analysis of anti-concentration of $\langle \Col_n(A),Y_n(A)\rangle$.
Recall that \textit{the L\'evy concentration function} $\cf(Z,t)$ of a real variable $Z$ is defined as
\[
\cf(Z,t):=\sup_{\lambda\in\R}\Prob\big\{|Z-\lambda|\leq t\big\},\quad t\geq 0.
\]
The relationship between the magnitude of $\cf\big(\sum_{i=1}^n a_i Z_i,t\big)$ for a linear combination
of random variables $\sum_{i=1}^n a_i Z_i$ and the structure of the coefficient vector $(a_1,\dots,a_n)$
has been studied in numerous works, starting from an inequality of Erdos--Littlewood--Offord \cite{LO1943,Erdos1945};
we refer, in particular, to works \cite{Rogozin1961,Kolmogorov1958,Kesten1969,Esseen1966}
as well as \cite{TV2009} and survey \cite{NV13} for a more recent account of \textit{the Littlewood--Offord theory}
and its applications to the matrix invertibility.

To characterize the structure of a coefficient vector in regard to anti-concentration,
the notion of the \textit{Essential Least Common Denominator} (LCD) has been introduced in \cite{RV08a}.
We quote a slightly modified definition from \cite{RV rect}:
\[
\LCD(a):=\inf\big\{\theta>0:\;\dist(\theta\,a,\Z^n)<\min(\gamma\lVert \theta\,a\rVert_2,\alpha\sqrt{n})\big\},\quad a\in\R^n.
\]
Here, $\alpha,\gamma$ are small positive constants.
The Littlewood--Offord--type inequality used in \cite{RV08a,RV rect} can be stated as follows.
If $Z_1,Z_2,\dots,Z_n$ are i.i.d real valued random variables with $\Prob\{|Z_i-\Exp\,Z_i|<\beta\}\leq 1-\beta$
for some $\beta>0$ then for any unit vector $a\in\R^n$,
\begin{equation}\label{3u32-4u4poiu}
\cf\bigg(\sum_{i=1}^n a_i Z_i,t\bigg)\leq Ct+\frac{C}{\LCD(a)}+2\exp(-c\,n),\quad t>0,
\end{equation}
where $C>0$ may only depend on $\beta,\gamma$ and $c>0$ only on $\alpha,\beta$
(see \cite{RV rect} for a proof).
Using an $\varepsilon$--net argument, the authors of \cite{RV08a} show that with probability exponentially close to one,
the random unit normal vector $Y_n(A)$ has an exponentially large $\LCD$.
This implies
\begin{align*}
\Prob\{|\langle \Col_n(A),Y_n(A)\rangle|\leq s\}
&\leq \Prob\big\{\LCD(Y_n(A))<\exp(c'n)\big\}+Cs
+2\exp(-c'\,n)\\
&\leq Cs+3\exp(-c''n),\quad s>0.
\end{align*}

The combination of all the ingredients now gives the final estimate
\[
\Prob\big\{s_{\min}(A)\leq t\,n^{-1/2}\big\}\leq \tilde C\,t+\tilde C\,\exp(-\hat c\,n),\quad t>0,
\]
matching, by the order of magnitude and up to the exponentially small additive term, the known asymptotics
of $s_{\min}$ of Gaussian random matrices \cite{Edelman,Szarek}.

\smallskip

In the remaining part of this section, we will consider some of the novel additions
to the methodology made over the past years. To avoid technical details as much as possible,
we will refer to compressible vectors as well as all related notions from the literature as \textit{almost sparse} vectors,
and incompressible vectors and their relatives as \textit{spread} vectors.

\smallskip

\textbf{Invertibility over almost sparse vectors.}
In the setting of \textit{dense} random matrices as described above,
the set of almost sparse vectors $\textrm{AlSp}_n$ can be treated by a simple $\varepsilon$--net argument
since anti-concentration estimates for $\lVert Az\rVert_2$ for an \textit{arbitrary} vector $z\in S^{n-1}$
are able to overpower the cardinality of the $\varepsilon$--net $\Net$ in $\textrm{AlSp}_n$.
In the case of sparse and certain models of
structured random matrices, such argument may not be sufficient since
the product $|\Net|\,\sup_{z\in S^{n-1}(\R)}
\Prob\{\lVert Az\rVert_2\leq s\}$ may become infinitely large even for small $s>0$.
We consider two [related] approaches to this problem from the literature.

A first one is based on further subdividing $\textrm{AlSp}_n$ into a few subsets $T_1,T_2,\dots$
according to the size
of set of vector's components of non-negligible magnitude,
and applying an $\varepsilon$--net argument within each of the subsets.
Anti-concentration estimates for $Az$ for vectors $z\in T_i$ then
compete with the cardinality of an $\varepsilon$--net on the set $T_i$ rather than on the entire collection $\textrm{AlSp}_n$,
which, for certain models, allows the proof to go through. We refer, in particular, to \cite[Section~4]{RZ2016}
and \cite[Section~3]{C16} for an application of this strategy to structured random matrices;
as well as \cite[Proposition~3.1]{Cook circ} dealing with adjacency matrices of random $d$--regular digraphs.

The second approach consists in identifying a class of non-random matrices $\mathcal{C}$ such that for every $M\in\mathcal {C}$
and every almost
sparse vector $z\in S^{n-1}$, $Mz$ has a non-negligible Euclidean norm, and then showing that with probability close to one, $A\in\mathcal{C}$.
As an example, consider a collection of matrices $M$
such that for every non-empty subset $I\subset[n]$ with $|I|\leq m$, there is a row $\Row_i(M)$ with
$|\supp\Row_i(M)\cap I|=1$. Then, it is not difficult to check that for every non-zero $m$--sparse vector $z$ one has
$Mz\neq 0$.
It can further be verified that a random matrix $A$ with i.i.d Bernoulli($p$) elements and $n^{-1}\textrm{polylog}(n)\leq p\leq c m^{-1}$
belongs to this class with probability tending to one with $n\to\infty$ \cite{BR16}.
The construction can be made robust to treat almost sparse vectors, and can be further elaborated to
deals with diagonal shifts of very sparse matrices \cite{BR16,LLTTY sing,RT circ}.

\smallskip

\textbf{Invertibility via distance.}
The relation \eqref{310984-198} discovered in \cite{RV08a} can be applied to any model of randomness.
However, this relation is not completely satisfactory when either (a) there are strong probabilistic dependencies
between $\Col_i(A)$ and $H_i(A)$ which make estimating $\Prob\big\{\dist(\Col_i(A),H_i(A))\leq t\}$ challenging,
or (b) invertibility over the almost sparse vectors cannot be treated with a desired precision
using approaches based on $\varepsilon$--net arguments or on conditioning on a particular structure
of the matrix.
Here, we consider some developments of the invertibility via distance argument made in the contexts of
$d$-regular random digraphs and smoothed analysis of the condition number.

Let $A_{n,d}$ be the adjacency matrix of a uniform random $d$-regular directed graph on $n$ vertices.
The regularity condition implies that for every $1\leq i\leq n$, $\Col_i(A_{n,d})$
is a function of $\{\Col_j(A_{n,d})\}_{j\neq i}$, creating issues with applying the original version of the argument
from \cite{RV08a}. 
In \cite{C14}, Cook proposed a modification of the argument based on considering distances between
the matrix columns and random subspaces of the form $H_{i_1,i_2,+}(A_{n,d}):=
\Span\{\Col_j(A_{n,d}),\,j\neq i_1,i_2;\;\Col_{i_1}(A_{n,d})
+\Col_{i_2}(A_{n,d})\}$, for $i_1\neq i_2$. That was later applied in \cite{Cook circ,LLTTY sing}.
Here, we quote \cite[Lemma~4.2]{LLTTY sing}:
denoting by $S(\rho,\delta)$ the collection of all unit vectors $x$ in $\Cp^n$ with $\inf\limits_{\lambda\in\Cp}
|\{i\leq n:\,|x_i-\lambda|>\rho\,n^{-1/2}\}|> \delta n$, one has
\begin{align*}
\Prob&\Big\{\inf_{x\in S(\rho,\delta)}\lVert (A_{n,d}-\,z\Id)x \rVert_2\leq t\,n^{-1/2}\Big\}\\
&\leq \frac{1}{\delta\,n^2}
\sum_{\substack{i_1,i_2\in [n],\\ i_1\neq i_2}}\Prob\big\{
\dist(\Col_{i_1}(A_{n,d}-\,z\Id),H_{i_1,i_2,+}(A_{n,d}-\,z\Id))\leq t/\rho\big\}.
\end{align*}
Conditioned on a realization of $\Col_j(A_{n,d})$, $j\neq i_1,i_2$ (hence, also $Y:=\Col_{i_1}(A_{n,d})
+\Col_{i_2}(A_{n,d})$), the support of the $i_1$--st column of $A_{n,d}$
is uniformly distributed on the collection of $d$-subsets $Q$ satisfying $\{j\leq n:\,Y_j=2\}
\subset Q\subset \supp Y$.
In the regime $d\to\infty$ as $n$ goes to infinity, this is ``enough randomness'' for a satisfactory
bound on $s_{\min}(A_{n,d}-\,z\Id)$ required by the Hermitization argument \cite{Cook circ,LLTTY sing}.

We remark here that another version of the argument for matrices with dependencies based on
evaluation of certain quadratic forms,
introduced in \cite{V sym}, has been used in a non-Hermitian setting in \cite{RV unit} to estimate the smallest
singular value of unitary and orthogonal perturbations of fixed matrices. We refer to \cite{RV unit} for details.

In \cite{Tikh shifted}, a variant of the invertibility via distance argument was developed to deal with
non-random shifts of matrices with continuous distributions.
A main observation of \cite{Tikh shifted} is that the distances $\dist(\Col_i(A),H_i(A))$, $1\leq i\leq n$,
are highly correlated, which allows for a more efficient analysis than the first moment method estimate \eqref{310984-198}.
The invertibility via distance is applied in \cite{Tikh shifted} to the entire sphere rather than the set of spread vectors.
As an illustration of the principle, we consider a simpler setting of centered random matrices
when the argument is still able to produce new results. 
Assuming $A$ is an $n\times n$ real random matrix with i.i.d entries of zero mean, unit variance, and the distribution
density bounded above by $\rho$, for every $t>0$ and $1\leq k\leq n$ one has
$
\Prob\big\{\exists\,I\subset[n]:\; |I|\geq k,\;\dist(\Col_i(A),H_i(A))\leq t\;\;\forall\,i\in I\big\}
\leq C_\rho\,t\,(n/k)^{5/11}
$, where $C_\rho>0$ may only depend on $\rho$ (see \cite[Prop.~3.8]{Tikh shifted}).
This, combined with the simple consequence of the \textit{negative second moment identity}
\[
s_{\min}(A)\geq \Big(\sum_{i=1}^n\dist(\Col_i(A),H_i(A))^{-2}\Big)^{-1/2},
\]
implies an estimate $\Prob\{s_{\min}(A)\leq t\,n^{-1/2}\}\leq C_\rho'\,t$, $t>0$,
which does not carry the $c^n$ additive term inevitable when an $\varepsilon$-net--based approach is used.
We refer to \cite{Tikh shifted} for the more involved setting of non-centered random matrices.

\smallskip

\textbf{Alternatives to the $\LCD$.} Functions of coefficient vectors different from
the Essential Least Common Denominator have been introduced in the literature
to deal with anti-concentration in the context of sparse and inhomogeneous random matrices as well as matrices with dependencies.
Here, we review some of them (for non-Hermitian models only).

The original notion of $\LCD$ is not applicable to the study of linear combinations of non-identically distributed variables:
in fact, given any vector $a\in S^{n-1}(\R)$ with an exponentially large $\LCD$,
one can easily construct mutually independent variables $Z_1,\dots,Z_n$ with $\cf(Z_i,1)\leq 1/2$, $1\leq i\leq n$,
and such that $\cf\big(\sum_{i=1}^n a_i Z_i,0\big)=\Omega(n^{-1/2})$.
Given a random vector $X$ in $\R^n$ and denoting by $\bar X$ the difference $X-X'$ (where $X'$ is an independent copy of $X$,
the \textit{Randomized Least Common Denominator} with respect to $X$ is defined as
\[
\RLCD^X(a):=\inf\big\{\theta>0:\;
\Exp\,\dist^2\big((\theta a_1 \bar X_1,\dots,\theta a_n \bar X_n),\Z^n\big)<
\min(\gamma\lVert \theta a\rVert_2^2,\alpha n)\big\},\; a\in\R^n.
\]
The notion was introduced in \cite{LTV} to deal with inhomogeneous random matrices with different
entries distributions. The small ball probability inequality \eqref{3u32-4u4poiu} from \cite{RV08a,RV rect} extends to the non-i.i.d setting
with the $\RLCD$ taking place of the original notion. We refer to \cite{LTV} for details.

Strong quantitative invertibility results for matrices with fixed rowsums
and adjacency matrices of $d$--regular digraphs obtained recently in \cite{Tran} and \cite{JSS reg}, respectively,
rely on a modification of the $\LCD$ which allows to treat linear combinations of Bernoulli
variables conditioned on their sum. Specifically, in \cite{Tran} the notion of the \textit{Combinatorial Least Common
Denominator} $\CLCD$ is defined as
\[
\CLCD(a):=\inf\big\{\theta>0:\;\dist(\theta (a_i-a_j)_{i< j},\Z^{n\choose 2})<\min(\gamma\lVert \theta (a_i-a_j)_{i< j}
\rVert_2,\alpha n)\big\},\;a\in\R^n,
\]
where $(a_i-a_j)_{i<j}$ denotes a vector in $\R^{n\choose 2}$ with $(i,j)$--th coordinate equal to $a_i-a_j$, $1\leq i<j\leq n$.
It is further shown that for the random vector $(Z_1,Z_2,\dots,Z_n)$ uniformly distributed on the collection of $0/1$ vectors
with exactly $n/2$ ones, an analog of the anti-concentration inequality \eqref{3u32-4u4poiu} holds,
with $\LCD$ replaced with $\CLCD$. A modification of the notion, called $\QCLCD$, was further considered in \cite{JSS reg}.
We refer to that paper for details.

Another functional --- \textit{the degree of unstructuredness} $\UD$ --- was introduced in \cite{LT2020}
to study invertibility of sparse Bernoulli random matrices. A main observation exploited in \cite{LT2020}
is that, for $p=o(1)$, linear combinations of i.i.d Bernoulli($p$) random variables
$\sum_{i=1}^n a_i Z_i$ are often more concentrated than corresponding linear combinations
of \textit{dependent} $0/1$ variables conditioned to sum to a fixed number of order $\Theta(pn)$.
In \cite{LT2020}, the argument proceeds by conditioning on the size of the support of a column of the matrix
and estimating anti-concentration of 
$\dist(\Col_i(A),H_i(A))=|\langle\Col_i(A),Y_i(A) \rangle|$ in terms of 
the degree of unstructuredness of the unit random normal $Y_i(A)$.
The definition of $\UD$
is technically involved, and we do not provide it here; see \cite{LT2020} for details.

\smallskip

\textbf{Average-case analysis of anti-concentration.}
The average-case study of Little\-wood--Offord--type inequlaities for linear combinations $\sum_{i=1}^n a_i Z_i$,
introduced in the random matrix context in \cite{Tikh2020}, was a crucial element in some recent advances
on quantitative invertibility of random discrete matrices \cite{Tikh2020,JSS2020a,JSS2020b},
which helped resolve some long standing problems in the combinatorial random matrix theory.
A main idea of \cite{Tikh2020} is, rather than attempting to obtain an explicit description of vectors $a$
such that $\sum_{i=1}^n a_i Z_i$ is strongly anti-concentrated, to consider the linear combination for
a \textit{randomly chosen} coefficient vector (with an appropriately defined notion of randomness).
This approach allowed to strengthen the invertibility results available through the use of the $\LCD$.
As an example, we consider a simplified version of the main technical result of \cite{Tikh2020}.
Let $\varepsilon\in(0,1/2)$, $M\geq 1$.
Then there exist $n_{0}=n_{0}(\varepsilon,M)$
depending on $\varepsilon,M$ and $L_{0}=
L_{0}(\varepsilon)>0$ depending \textit{only} on $\varepsilon$
(and not on $M$)
with the following property. Take $n\geq n_{0}$, $1\leq N\leq (1/2+\varepsilon)^{-n}$,
and let $\mathcal A:=(\{-2N,\dots,-N-1\}\cup\{N+1,\dots,2N\})^n$.
Assume that a random vector $a=(a_1,\dots,a_n)$ is uniformly distributed on $\mathcal A$.
Then
\[
\Prob_a\big\{\cf_Z \big(a_1 Z_1+\dots+a_n Z_n,\sqrt{n}\big)
> L_{0} N^{-1}
\big\}\leq e^{-M\,n}.
\]
Here, $\cf_Z(\cdot,\cdot)$ denotes the L\'evy concentration function with respect to the randomness of $(Z_1,\dots,Z_n)$,
a vector with independent $\pm 1$ components. 
The main point of the statement is that the parameter $L_0$ controlling the anti-concentration of the linear combination,
does not depend on $M$, i.e the proportion of the coefficient vectors in $\mathcal A$
such that anti-concentration of $a_1 Z_1+\dots+a_n Z_n$ is weak, becomes \textit{superexponentially small} in $n$
as $n\to\infty$. 

\smallskip

\textbf{Matrices with heavy entries.} For invertibility of (dense) random matrices with independent entries assuming only
finite second moments, we refer to \cite{RT15,Livshyts,LTV}.

\section{Open problems}\label{s:open}

We conclude this survey with a selection of open research problems. 

\textbf{Refined smoothed analysis of invertibility.}
Recall that a standard model in the setting of the smoothed analysis of the condition number
is of the form $A+M$, where $A$ is an $n\times n$ random matrix with i.i.d entries,
and $M$ is a non-random shift.
\begin{problem}[Shift-independent estimates for matrices with continuous distributions]
Let $\xi$ be a real random variable of zero mean, unit variance, and bounded distribution density.
Let $A$ be an $n\times n$
matrix with i.i.d entries equidistributed with $\xi$.
It is true that for every non-random matrix $M$,
\[
\Prob\big\{s_{\min}(A+M)\leq t\,n^{-1/2}\big\}\leq Ct,\quad t>0,
\]
where $C>0$ may only depend on the c.d.f of $\xi$ (and not on $n$)?
\end{problem}
For partial results on the above problem, see \cite{SST06,Tikh shifted}.

\smallskip

\begin{problem}[Optimal dependence of $s_{\min}(A+M)$ on the norm of the shift in the discrete setting]
Let $A$ be an $n\times n$
matrix with i.i.d $\pm1$ entries, and let $T,t>0$ be parameters.
For any $\varepsilon,L>0$, estimate
$
\sup_{M:\,\lVert M\rVert\leq T}\Prob\big\{s_{\min}(A+M)\leq t\big\}
$
up to a multiplicative error $O(n^{\varepsilon})$ and an additive error $O(n^{-L})$, that is,
find an explicit function $f(n,T,t)$ such that
\[
n^{-\varepsilon}\,f(n,T,t)-C\,n^{-L}\leq
\sup_{M:\,\lVert M\rVert\leq T}\Prob\big\{s_{\min}(A+M)\leq t\big\}\leq n^{\varepsilon}\,f(n,T,t)+C\,n^{-L},
\]
where $C>0$ may only depend on $\varepsilon$ and $L$.
\end{problem}
For best known partial results on the above problem, see \cite{TV10c,JSS2020 smooth}.

\smallskip

\begin{problem}[Dependence of $s_{\min}(A+M)$ on $M$ in the Gaussian setting]
Let $A$ be an $n\times n$ matrix with i.i.d standard real Gaussian entries. Find an estimate
on $\Exp\,s_{\min}(A+M)$ in terms of the singular spectrum of $M$.
\end{problem}
One can assume in the above problem
that $M$ is a diagonal matrix with $i$-th diagonal element $s_i(M)$, $1\leq i\leq n$.
Note that $A$ may either improve or degrade invertibility of $M$.

\medskip

\textbf{Invertibility and spectrum of very sparse matrices.} Here, we consider the problem
of identifying the limiting spectral distribution for non-Hermitian matrices with \textit{constant} average
number of non-zero elements in a row/column.

\begin{problem}[{The oriented Kesten--McKay law; see \cite[Section~7]{BC2011}}]
Let $d\geq 3$. For each $n$, let $A_{n,d}$ be the adjacency matrix of a uniform random $d$--regular directed
graph on $n$ vertices. Prove that the sequence of spectral distributions $\big(\mu_{A_{n,d}}\big)_{n=1}^\infty$
converges weakly to the probability measure on $\Cp$ with the density function
\[
\rho_d(z):=\frac{1}{\pi}\frac{d^2(d-1)}{(d^2-|z|^2)^2}\,
\textbf{1}_{\{|z|<\sqrt{d}\}}.
\]
\end{problem}
Assuming the standard Hermitization approach to the above problem,
the following is the crucial (perhaps the main) step of the argument:
\begin{problem}
Let $d\geq 3$ and let $(A_{n,d})_{n=1}^\infty$ be as above. Prove that for almost every $z\in\Cp$ and every $\varepsilon>0$,
\[
\lim_{n\to\infty}\Prob\big\{s_{\min}(A_{n,d}-z\,\Id)\leq \exp(-\varepsilon\,n)\big\}=0.
\]
\end{problem}

\smallskip

\begin{problem}[Spectrum of directed Erdos--Renyi graphs of constant average degree]
Let $\alpha>0$. For each $n\geq \alpha$, let $A_n$ be an $n\times n$ random matrix with i.i.d Bernoulli($\alpha/n$) entries.
Does a sequence of spectral distributions $\big(\mu_{A_n}\big)$ converge weakly to a non-random probability measure?
\end{problem}
As in the case of regular digraphs, assuming the Hermitization argument, the following problem constitutes
an important step to understanding asymptotics of the spectrum:
\begin{problem}
For each $n\geq \alpha$, let $A_n$ be an $n\times n$ random matrix with i.i.d Bernoulli($\alpha/n$) entries.
Is it true that for almost every $z\in\Cp$ and every $\varepsilon>0$,
\[
\lim_{n\to\infty}\Prob\big\{s_{\min}(A_n-z\,\Id)\leq \exp(-\varepsilon\,n)\big\}=0?
\]
\end{problem}

\medskip

\textbf{Invertibility and spectrum of structured random matrices.}
The spectrum of structured random matrices in the absence of 
expansion-like properties (such as \textit{broad connectivity} \cite{RZ2016} or \textit{robust irreducibility} \cite{Cooketal1,Cooketal2})
is not well understood as of now. In particular, a full description of the class
of inhomogeneous matrices with independent entries with spectral convergence to the circular law
seems to be out of reach of modern methods.
\begin{problem}
Give a complete description of sequences of standard deviation profiles $(U_n)_{n=1}^\infty$ satisfying the following condition:
assuming that $\xi$ is any random variable of zero mean and unit variance, and that for each $n$,
$M_n$ is an $n\times n$ matrix with i.i.d entries equidistributed with $\xi$, the sequence of spectral distributions
$\big(\mu_{U_n\odot M_n}\big)$ converges weakly in probability to the uniform measure on the unit disc of $\Cp$.
\end{problem}
A natural class of profiles considered, in particular, in \cite{Cooketal1,Cooketal2}, are \textit{doubly stochastic} profiles.
One may expect that those profile sequences, under some weak assumption on the magnitude of the maximal entry,
should be sufficient for the circular law to hold:
\begin{problem}
Assume that for each $n$, the standard deviation profile $U_n$ satisfies
\[
\sum_{i=1}^n (U_n)_{ij}^2=\sum_{i=1}^n (U_n)_{ji}^2=1,\quad 1\leq j\leq n,
\]
and that for some $\varepsilon>0$, $\limsup_n\max_{ij}((U_n)_{ij}n^\varepsilon)=0$.
Is it true that, with $M_n$ as in the above problem, the sequence
$\big(\mu_{U_n\odot M_n}\big)$ converges weakly in probability to the uniform measure on the unit disc of $\Cp$?
\end{problem}
Note that the above setting allows sparse matrices (cf. \cite[Theorem~2.4]{Cooketal1}).
Solution to the above problem, if approached
with Girko's Hermitization procedure, requires satisfactory bounds on the smallest singular values of
$U_n\odot M_n-z\,\Id$.




\end{document}